\documentclass{article} 
\usepackage[all]{xy}
\usepackage{amsmath}
\usepackage{amssymb}
\usepackage{amsthm}

\newtheorem{thm}{Theorem}[section]
\newtheorem{lem}[thm]{Lemma}
\newtheorem{cor}[thm]{Corollary}
\newtheorem{prop}[thm]{Proposition}

\newtheorem{defn}{Definition}

\newtheorem{rem}[thm]{Remark}
\newtheorem{ex}[thm]{Example}

\def\C{{\mathbb C}}
\def\D{{\mathbb D}}

\def\P{{\mathbb P}}

\def\cE{{\cal E}}
\def\cF{{\cal F}}

\def\cL{{\cal L}}
\def\cM{{\cal M}}

\def\cO{{\cal O}}

\def\iso{\cong}

\def\lra{\longrightarrow}

\def\ra{\rightarrow}

\def\tensor{\otimes}
\def\operatorname#1{\mathop{\rm #1}\nolimits}

\def\Chow{\operatorname{Chow}}

\def\deg{\operatorname{deg}}

\def\Chow{\operatorname{Chow}}
\def\RatCurves{\operatorname{RatCurves}}

\title{Holomorphic Engel Structures}
\author{ Francisco Presas \and Luis Eduardo Sol\'a Conde }

\begin{document}

\maketitle

\noindent {\bf Abstract:}
Recently there has been renewed interest among differential geometers in the theory of Engel structures. We introduce holomorphic analogues of these structures, and pose the problem of classifying projective manifolds admitting them. Besides providing the basic properties of these varieties, we present two series of examples and characterize them by certain positivity conditions on the Engel structure.

% The theory of been a growing interest in Engel structures within the have beenIn this paper we introduce the concept of holomorphic Engel structure, the complex analogue of the classic real Engel structuresThe study of Engel structures over complex $4$-manifolds is started. Basic properties and examples are provided. Moreover, a complete classification result is outlined.

\tableofcontents

%%%%%%%%%%%%%%%%%%%%%%%%%%%%%%%%%%%%%%%%%%%%%%%%%%%%%%%%%%%%%
\section{Introduction.}\label{sec:intro}

A class of distributions of a manifold is
called open if it is open as a set inside the parameter space of distributions of
the manifold.
For instance the class of contact distributions is open, in the real
and holomorphic cases. Moreover, the class is called {\it topologically
stable} if it is open and it admits a unique local model in the neighborhood
of each point. Again contact distributions are topologically stable
since any of them admits a holomorphic chart $(z, x_1, y_1, \ldots
x_n, y_n)$ at a given point in which it becomes $\ker (dz-\sum x_idy_i)$.
In the real case, it has been proved (cf. \cite{Mo1, Mo2}) that the list of topologically
stable distributions is given by:
\begin{itemize}
\item line fields,
\item contact structures in odd dimensional manifolds,
\item even contact structures in even dimensional manifolds,
\item Engel structures in $4$-dimensional manifolds.
\end{itemize}
The definitions of the last two ones will be provided later. We must remark that in the holomorphic case a similar proof provides the same list.
The name of "Engel structure" comes from the fact that H. Engel gave the first proof of the topological stability of
this kind of distributions in the 19th century, basically by showing that they follow a canonical local model, i.e a kind of Darboux theorem. By completeness, in Subsection \ref{sub:darboux}, we offer the proof of the equivalent result in the holomorphic case, which guarantees the topological stability of the holomorphic Engel structures.

An Engel structure on a (real/complex) variety $X$ is a
$2$-dimensional distribution $D\subset T_X$
satisfying that $[D, D]= \cE$ is a rank $3$ distribution and $[\cE, \cE]=T_X$. It is simple to check that there is a
one dimensional subbundle $L \subset \cE$ that is invariant for the
bracket operation, i.e. $[L, \cE] \subset \cE$. Moreover, we have $L
\subset D$. This provides a flag
$$ L \subset D \subset \cE \subset T_X$$
of vector bundles inside the tangent bundle. In the real orientable case, this
immediately implies that the tangent bundle is parallelizable. A
long standing conjecture in differential geometry has been that the
converse is also true: any parallelizable $4$-fold admits an Engel
structure. This has been recently proved by A. Vogel \cite{V},
lighting the interest for (real) Engel structures.

We want to start in this note the study of holomorphic Engel
structures on projective varieties. The main result we show is
\begin{thm}\label{main-thm}
Let $X$ be a projective variety admitting an Engel structure $L\subset D\subset\cE$. Then either:
\begin{itemize}
\item $L^{-1}$ is not pseudoeffective and then $X$ is the Cartan prolongation of a contact $3$-fold, or
\item $(D/L){-1}$ is not pseudoeffective. In this case, if we assume moreover that $D\cong L\oplus D/L$, then $X$ is a Lorentzian tube.
\end{itemize}
Moreover the two classes have a unique common element which is the universal
family of lines contained in a quadric hypersurface in $\P^4$.
\end{thm}
The definitions of Cartan prolongation and Lorentzian tube will be given in section \ref{sec:ex}. It is enough for
now to understand that the Engel structures described in the previous list are classified by $3$-dimensional contact
structures and $3$-dimensional conformal structures. Remark that
$3$-dimensional contact structures are completely classified, since
the abundance conjecture \cite{KM} is true in dimension $3$ and then
by \cite{KPSW} $X$ is either a projectivized tangent bundle or a Fano
$3$-fold. On the
other hand holomorphic conformal structures also have been
classified in a recent article \cite{JR}. So, the Theorem
\ref{main-thm} provides a very concrete classification result.

The structure of this note is as follows. In Section
\ref{sec:engel} we will provide the basic definitions and
properties we will play with, also proving the local stability result mentioned as ``Darboux Theorem''. In Section \ref{sec:ex} the two classical examples of Theorem \ref{main-thm} will be explained. Finally we will
provide the proof of Theorem \ref{main-thm} in Section \ref{sec:thm}.

\medskip

\noindent {\bf Thanks and acknowledgments.}
Very special thanks to Dan Fox, who pointed out to us the second construction of Engel structures due to Cartan and not known for most of the people working in the area. We want to thank Ignacio Sols and Tom\'as Luis G\'omez for fruitful discussions.

%%%%%%%%%%%%%%%%%%%%%%%%%%%%%%%%%%%%%%%%%%%%%%%%%%%%%%%%%%%%%
\section{Basic properties}\label{sec:engel}

\subsection{Global properties.}

Along this paper $X$ will denote a smooth projective complex variety
of dimension $2n$, and $T_X$ its tangent bundle. Given any subsheaf
$E\subset T_X$, the composition of the usual Lie-bracket of $E$ with
the projection onto $T_X/E$
$$
E\tensor E\to T_X/E$$
is an $\cO_X$-linear morphism, that we will
call { \it O'Neill tensor} of $E$.

\begin{defn}\label{defn:even}{\rm
    Let $\cE\subset T_X$ denote a codimension $1$ subbundle, defined
    as the kernel of a surjective morphism:
    $$
    T_X\stackrel{\theta}{\lra}L'\ra 0.
    $$
    Locally we can compute the exterior differential $d\theta$, but
    these data do not glue together to define a $2$-form. Nevertheless
    a direct computation shows that $\theta\wedge (d\theta)^{n-1}$ is a well
    defined $(2n-1)$-form with values in $(L')^{n}$.  We say that $\cE$ defines
    an {\it even contact structure} if $\theta\wedge (d\theta)^{n-1}\in
    H^0(\Omega_X^{2n-1}\tensor (L')^{n})$ is everywhere non-zero.}
\end{defn}

\begin{rem}\label{rem:even}{\rm
    Furthermore $d\theta$ defines a global section of
    $\bigwedge^2\cE^\vee\tensor L'$, that we denote by the same
    symbol. In fact $d\theta$ coincides with the O'Neill tensor via
    the Cartan formula
    $$
    d\theta(V_1,V_2)=\theta(V_1)V_2-\theta(V_2)V_1-\theta[V_1,V_2].
    $$}
\end{rem}

\begin{rem}\label{rem:L}
  {\rm By the non-vanishing condition, the kernel of the skew-symmetric morphism
    $d\theta:\cE\to\cE^\vee\tensor L'$ is a line subbundle of $\cE$,
    that we denote by $L$. We call it the {\it kernel of the
      even contact structure}. With this notation $d\theta$ provides
    an everywhere nondegenerate skew-symmetric morphism
    $\bigwedge^2(\cE/L)\to L'$. In particular, if $\dim X=4$, $\det(\cE)=L'\tensor
    L^{-1}$ and $\omega_X^{-1}=(L')^{2}\tensor L^{-1}$.}
\end{rem}

\begin{defn}\label{def:engel}{\rm
    An {\it Engel structure} is a rank $2$
    vector subbundle $D\subset T_X$ on a $4$-dimensional complex manifold $X$ satisfying that $\cE=[D,D]$
is an even contact structure. Here $[,]$ denotes the usual Lie bracket.}
\end{defn}

%Recall that $\cE$ is not part of the definition since it is built
%out from $D$.

\begin{lem}\label{lem:engel}
  {\rm An Engel structure provides a filtration of $T_X$:
    $$
    L\subset D\subset\cE\subset T_X.
    $$
    Moreover the O'Neill tensor of $D$ induces an isomorphism
    $L\tensor D/L\iso\cE/D$.}
\end{lem}

\begin{proof}
  The statement is local, hence we assume that $D$ is generated by two
  vector fields $V_1$ and $V_2$, where $\{V_1,V_2,V_3:=[V_1,V_2]\}$
  generates $\cE$. We may assume that $[V_2,V_3]\notin\cE$.  A vector
  field $V=aV_1+bV_2+cV_3$ belongs to $L$ iff $[V_i,V]\in\cE$,
  $i=1,2,3$.  Taking $i=2$ we obtain $c=0$, therefore $L\subset D$.

  We may now assume that $L$ is generated by $V_1$, and the class of
  $V_2$ generates $D/L$. Then the class of $V_3=[V_1,V_2]$ generates
  $\cE/D$.
\end{proof}

We will usually identify an Engel structure with its associated filtration.
The elements of this filtration fit into
the following commutative diagram with exact rows and columns
\begin{equation}
\xymatrix{
  & 0 \ar[d] & 0 \ar[d] &  &  \\
  0 \ar[r] & L \ar@{=}[r] \ar[d]& L \ar[d] \ar[r] & 0 \ar[d]
  & \\
  0 \ar[r] & D \ar[r]\ar[d] & \cE \ar[r]\ar[d]
  & L\tensor D/L \ar[r] \ar@{=}[d] & 0 \\
  0 \ar[r] & D/L \ar[r]\ar[d] & \cE/L \ar[r]\ar[d]
  & L\tensor D/L \ar[r] \ar[d] & 0 \\
  & 0 & 0 & 0 &  \\
}\label{diag1}
\end{equation}

\begin{rem}
  {\rm It does not seem true in general that an even contact structure determines
  an Engel structure. Moreover, there are examples of even contact
  structures not supporting an Engel structure at all (see Corollary
  \ref{cor:non-engel})}.
\end{rem}

\begin{lem}\label{lem:canon}
  With the above notation, the anticanonical line bundle
  $\omega_X^{-1}$ is isomorphic to $L^3\tensor (D/L)^4$.
\end{lem}
\begin{proof}
  It follows directly from Remark \ref{rem:L} and the commutative
  diagram above:
  $$
  \omega_X^{-1}=L\tensor\det(\cE/L)\tensor
  L'\stackrel{\ref{rem:L}}{=}L\tensor\det(\cE/L)^2=L\tensor\big(L\tensor
  (D/L)^2\big)^2=L^3\tensor (D/L)^4
  $$
\end{proof}

\subsection{Local properties.} \label{sub:darboux}
Now, we will prove the existence of a local canonical form for any Engel structure.

\begin{thm}[Darboux Theorem] \label{engel_local}
Let $D$ be an Engel structure on $X$ and let $x$ be a point on $X$. Then,
there exists a neighborhood $U$ of $x$ in $X$ and a
chart $\phi: V \subset \C^4 \to U$ such that $\phi^* D= D_0$,
where $D_0$ is the Engel structure in $\C^4$ given by the
vanishing conditions
\begin{eqnarray}
dx-ydz & = & 0 \label{eq:En1}\\
dy -wdx & = & 0 \label{eq:En2}
\end{eqnarray}
\end{thm}

We need a couple of holomorphic contact geometry results before addressing the proof of this Theorem.

\begin{lem}[Holomorphic Gray stability lemma] \label{le:Gray}
 Let $X$ be a smooth variety, let $U\subset \C$ be a neighborhood of $0$, let $x$ be a point in  $X$ and assume that there is a family of contact forms
$$
\cF_t \to T_X \stackrel{\theta_t}{\to} \cO, $$ with $t\in U \subset
\C$. Then there is a neighborhood $V$ of $x\in X$ and family of holomorphic flows $\phi_t: V \to X$ such that
$\phi_t^* \cF_t = \cF_0$, for $|t|$ small enough.
\end{lem}
The flow is not well-defined for all $t\in U$, but the
holomorphic vector fields generating it are. The point is that since
we do not assume compactness the fields do not always integrate.
However, this will not be a problem in the applications. We are
assuming that $L= \cO$ to make the proof simpler; this case is
enough for our purposes, though the proof can be probably extended to other situations.
\begin{proof}
 The condition is equivalent to the following one
\begin{equation} \label{gray}
 \phi_t^* \theta_t = f_t \theta_0,
\end{equation}
for some $f_t \in H^0(\cO^*)$. The flow $\phi_t$ is uniquely
represented by a family of holomorphic vector fields $V_t$.
Differentiating in the previous equation we obtain
$$
\phi_t^* \left( \frac{d}{dt} \theta_t + \cL_{V_t} \theta_t \right) =
\frac {df_t}{dt} \theta_0,
$$
which can be rewritten as
\begin{equation} \label{eq:auxi}
\frac{d}{dt} \theta_t + \cL_{V_t} \theta_t = h_t \theta_t,
\end{equation}
with $h_t= (\frac{1}{f_t} \frac{d}{dt} f_t) \circ \phi_t^{-1}$.
There is a well-defined non-zero holomorphic vector field associated
to any contact form, whenever $L= \cO$, that is called the Reeb field
$R$ and it is defined by the following two equations
\begin{eqnarray*}
 i_R d\theta & = & 0 \\
i_R \theta & = & 1.
\end{eqnarray*}
Recall that $\ker d\theta$ is a $1$-dimensional space, so the first
condition determines a $1$-dimensional foliation and the second one
is just a normalization. The Reeb fields associated to our family of
forms $\theta_t$ will be denoted by $R_t$. We impose $h_t= i_{R_t} \frac{d}{dt}
\theta_t$. Note that $\cF \bigoplus R_t = T_X$. It is simple to
check that the previous assignation makes the equation (\ref{eq:auxi})
true for any vector in the direction $R_t$. So, we have to check it
just for vectors $Y\in \cF_t$. Use Cartan's formula
$$ \cL_Y= d i_X + i_X d, $$
and assume that we check it against $Y\in \cF_t$, then we have
$$
\frac{d}{dt} \theta_t + d i_{V_t} \theta_t + i_{V_t} d\theta_t = 0.
$$ We just look for a vector $V_t \in \cF_t$ and this further reduces the
expression to
$$
\frac{d}{dt} \theta_t + i_{V_t} d\theta_t = 0. $$ Now, recall that
$d\theta_t$ is non degenerate on $\cF_t$, so this equation has a
unique solution $V_t$. This ends the proof.
\end{proof}

\begin{cor}  \label{cor:Darboux}
 Let $(X, \cF)$ be a contact manifold. For any point $x\in X$, there exists a neighborhood $U$
 and a chart $\phi: V \subset \C^{2n+1} \to U$ such that $\phi^* \cF= \cF_0$, where $\cF_0$ is the
 standard contact form in $\C^{2n+1}$ defined by the vanishing condition
\begin{equation*}
dx-\Sigma_{j=1}^n y_jdz_j = 0.
\end{equation*}
\end{cor}
\begin{proof}
Take any chart $\phi: V \subset \C^{2n+1} \to U$ centered at $x$.
Choose an element of $A \in GL(2n+1, \C)$ and construct $\Phi= A
\circ \phi$ in such a way that $\Phi^* \theta(x) = \theta_0(0)$ and
$\Phi^* d\theta(x) = d\theta_0(0)$. This can be done just by
choosing symplectic basis at $\cF(x)$ and $\cF_0(0)$ and choosing
$A$ in order to map the first one into the second one. Now, to be contact is
an open condition so $\theta_t = (1-t)\theta_0 + t\Phi^* (\theta)$
is a family of contact structures for some small neighborhood $W$ of
$0\in \C^{2n+1}$. We are in the hypothesis of Lemma \ref{le:Gray}
and therefore there exists a flow, for $|t|$ small enough, $\Psi_t: W \to
W$ such that $\Psi_t^* \ker (\theta_t) = \cF_0$. In fact, probably
shrinking $W$ again, the flow is well defined for $t=1$. We conclude
that $\Psi_1 \circ \Phi$ is the needed map.
\end{proof}

Now we complete the proof of Theorem \ref{engel_local}.

\begin{proof}
We will construct a special chart $(\phi, U)$ at $p\in X$ with local coordinates $(x,y,z,w)$. The first condition imposed to the chart
will be to send the line bundle $L$ to the derivative
$\frac{\partial}{\partial w}$ of a function $w:U \to \C$.  The second and final one will be to choose local coordinates $(x,y,z)$ in a hypersurface $w= \epsilon$ in such a way that the coordinates $(x,y,z,w)$ are provided by the flow along the field $\frac{\partial}{\partial w}$.

Now, in these special coordinates the bundle $\cE$ is $w$-invariant so projects to
the quotient $U/L$, which has coordinates $(x,y,z)$, as a bundle
$\hat{\cE}$. This bundle defines a contact structure on the quotient. Denote the projection
by $\pi: U \to U/L$. We apply Corollary \ref{cor:Darboux} to $U/L$ to obtain a map $\hat{\phi}: V/L \to
V/L$, where $V\subset U$, satisfying that $\hat{\phi}(\hat{\cE})=
\hat{\cE_0}$. The last one is the standard contact structure in
$\C^3$ and is defined by the condition (\ref{eq:En1}). Denote by $\phi$ the fiberwise-constant lift of $\hat{\phi}$
to $V$. Let $q$ be a point on $V$. Recall that $D_q \supset L_q$ can be
projected to $V/L$ defining an element in $\P(\hat{\cE}_{\pi(q)})$. This provides
a map $\Psi: V \to W \subset \P(\hat{\cE})$, where $W$ is an open subset.
The condition of being Engel is completely equivalent to the biholomorphicity of this map and this implies that
the Engel structure is determined by an equation of the form
$$ dy- gdx=0, $$
with $g$ a certain holomorphic function satisfying that $\frac{\partial g}{\partial w} \neq 0$.
By another change of coordinates, provided by the Inverse Function Theorem, the Engel structure gets defined by condition (\ref{eq:En2}).

\end{proof}

%%%%%%%%%%%%%%%%%%%%%%%%%%%%%%%%%%%%%%%%%%%%%%%%%%%%%%%%%%%%%
\section{Examples}\label{sec:ex}

We have seen in the previous section, thanks to the Darboux Theorem, that every local Engel structures is defined, up to holomorphic change of coordinates, by the vanishing of the equations
\begin{eqnarray*}
dx-ydz & = & 0 \\
dy -wdx & = & 0.
\end{eqnarray*}
In this section we will present two compact examples of Engel structures. They are holomorphic analogues of two real constructions already known to Cartan. The other two known examples of real Engel structures, a construction over mapping tori \cite{Ge} and a general construction on parallelizable $4$-folds \cite{V}, cannot be adapted to the holomorphic case.
%\begin{ex}\label{ex:local}
%  {\rm Take the following two equations
%  \begin{eqnarray*}
%dx-ydz & = & 0 \\
%dy -wdx & = & 0.
%\end{eqnarray*}
%They determine an Engel structure in $\C^4$ as a trivial computation
%shows. In fact, we show later, that locally an Engel structure has
%always this form.}
%\end{ex}

\begin{ex}[Cartan prolongation of a contact structure]\label{ex:main}
  {\rm Let $Z$ be a projective smooth threefold admitting a contact
    structure $F\subset T_Z$ appearing as the kernel of a morphism
    $\varphi:T_Z\to L'_Z$. Let $X$ be the projectivization of
    $F^\vee$, and denote by $\pi$ the natural projection onto $Z$.
    Composing the pull-back of $\varphi$ with the differential of $\pi$ we obtain a
    twisted $1$-form
    $$
    \theta:=\pi^*\varphi\circ d\pi:T_X\to L':=\pi^*L'_Z.
    $$
    Being $F$ a contact distribution, it follows that
    $\theta$ defines an even contact structure $d\pi^{-1}(\pi^*F)$ on
    $X$, that we denote by $\cE$, and in fact the converse is also true, as one can easily check locally: if the contact structure on $Z$ is given locally by the $1$-form $\theta_Z$, then $\theta\wedge d\theta=\pi^*(\theta_Z\wedge d\theta_Z)$. Note also that, by construction, $\ker(\cE)=T_{X|Z}$.

    Furthermore, we claim that it admits an Engel structure. In order
    to see that, consider the relative Euler sequence associated with
    the tautological line bundle $\cO_X(1)$, twisted with $\cO_X(-1)$:
    $$
    \xymatrix{0\ar[r] & \cO_X(-1)\ar[r] & \pi^*F \ar[r]
      &T_{X|Z}(-1)\ar[r] &0.}
    $$
    Denote by $D$ the kernel of the composition $\cE\to\pi^*F\to
    T_{X|Z}(-1)$. We get the following commutative diagram with exact
    rows and columns:
    $$
    \xymatrix{
      & 0 \ar[d] & 0 \ar[d] &  &  \\
      0 \ar[r] & L \ar@{=}[r] \ar[d]& T_{X|Z} \ar[d] \ar[r] & 0 \ar[d]
      & \\
      0 \ar[r] & D \ar[r]\ar[d] & \cE \ar[r]\ar[d]
      & T_{X|Z}(-1) \ar[r] \ar@{=}[d] & 0 \\
      0 \ar[r] & \cO_X(-1) \ar[r]\ar[d] & \pi^*F=\cE/L \ar[r]\ar[d]
      &  T_{X|Z}(-1) \ar[r] \ar[d] & 0 \\
      & 0 & 0 & 0 & \\}
    $$}
\end{ex}

Note that, in a more general setting, a regular fibration by curves over a contact manifold supports an even contact structures whose kernel is the relative tangent bundle of the fibration. Hence it makes sense to ask whether these even contact structures are associated to Engel structures.
The next proposition and its corollary show that this is only possible if the fibers have genus $0$.

\begin{prop}\label{prop:fibration}
  Let $X$ be projective smooth $4$-fold admitting an Engel structure
  $L\subset D\subset\cE$. Assume that all the leaves of $L$ are
  compact.  Then the genus of the leaves is zero.
\end{prop}

\begin{proof}
  Let $C$ be a general compact leaf of the foliation $L$. Since $X$ is
  a fibration locally around $C$, then the normal bundle of $C$ in $X$
  is trivial. This implies that $\omega^{-1}_X|_C=T_C=L|_C$, hence by
  Lemma \ref{lem:canon} the divisor
  $L'|_C=\det(\cE/L)|_C=\big(L\tensor (D/L)^2\big)|_C$ has degree
  zero. On the other side the exact sequence
  $$
  \xymatrix{0\ar[r]&(\cE/L)|_C\ar[r]&(T_X/L)|_C\iso\cO_C^{\oplus
      3}\ar[r]&L'|_C=(T_X/\cE)|_C\ar[r]&0}
  $$
  tells us that $L'|_C$ is globally generated, therefore trivial.
  It follows that $(\cE/L)|_C$ is trivial, too.

  Next we consider the exact sequence appearing in the last row of
  diagram \ref{diag1}, restricted to $C$:
  $$
  \xymatrix{0\ar[r]&(D/L)|_C\ar[r]&(\cE/L)|_C\iso\cO_C^{\oplus
      2}\ar[r]&(L\tensor D/L)|_C\ar[r]&0.}
  $$
  Note that the triviality of $\big(L\tensor (D/L)^2\big)|_C$
  implies that $\deg(L\tensor D/L)|_C=1-g$. But a quotient of a
  trivial bundle cannot have negative degree, hence $g\leq 1$.

  Finally if $g=1$ the divisors $L$ and $D/L$ are trivial along the
  leaves of $L$. Being $Z$ the subscheme of $\Chow(X)$ parameterizing these leaves and
  $\pi:X\to Z$ the projection, $D/L$ is the pull-back of a line
  subbundle $L_Z\subset T_Z$. Hence $D$ is the inverse image of $\pi^*
  L_Z$ by $d\pi$, and so $D$ is trivial along the fibers of $\pi$,
  contradicting its non integrability.

\end{proof}

\begin{cor} \label{cor:non-engel}
  Let $Z$ be a $3$-dimensional contact variety. Let $X \to Z$ be a
  regular fibration by curves. If the genus of the curves is not zero,
  then the even contact structure defined on $X$ by pull-back does not
  support an Engel structure.
\end{cor}

\begin{rem}\label{rem:non-engel}
{\rm In the setting of Proposition \ref{prop:fibration}, it follows that $L$ has positive degree on its leaves. In particular Proposition \ref{prop:cartan} below will imply that $X$ is in fact the Cartan prolongation of a contact manifold. }
\end{rem}

The following example is based on an analogous construction in (real)
differential geometry due to Cartan, where an Engel structure is built
upon a Lorentzian structure on a threefold. Cartan did not show that the structure was Engel but studied its properties.

We recall some definitions from the literature (see for instance \cite{JR}).
\begin{defn} \label{def:conformal}{\rm
A {\it holomorphic conformal structure} $\lambda$ on a holomorphic bundle
  $E$ over a complex manifold $Y$ is a non-degenerate section $\lambda\in
  H^0(Y, S^2 E^* \otimes M)$, where $M$ is an arbitrary line bundle.
  A {\it holomorphic conformal structure on a complex manifold} is a
  holomorphic conformal structure on the tangent bundle of the
  manifold.}
\end{defn}

\begin{defn}[Tubes \cite{GS}]\label{def:tube}
  {\rm
A {\it tube} on a complex manifold $X$ is the zero set of the
  defining section of a holomorphic conformal structure over a rank
  $r$ vector bundle $E$. It has a natural structure of
  $Q_{r-2}$-fibration over $X$
A {\it Lorentzian tube} is a tube on the tangent bundle of the
  manifold.}
\end{defn}

\begin{ex}[Lorentzian tubes]\label{ex:second}
  {\rm Let $Y$ be a smooth threefold equipped with a Lorentzian tube
    provided by a section $\lambda\in H^0(Y, S^2\Omega_Y\otimes M)$. Equivalently, we have a smooth conic fibration
    $\rho:X:=\mbox{zeroes}(\lambda)\subset\P(\Omega_Y)\to Y$. Consider
    the rank $2$ subbundle $D\subset T_X$ defined in the following
    way: being $i:X\hookrightarrow\P(\Omega_Y)$ the natural inclusion
    and denoting $L:=\big(i^*\cO_{\P(\Omega_Y)}(-1)\big)$ consider the
    diagram
    $$
    \xymatrix{&&T_X\ar[d]^{d\rho}&&\\0\ar[r]&L\ar[r]&\rho^*T_Y\ar[r]&i^*T_{\P(\Omega_Y),Y}(-1)\ar[r]&0,}
    $$
    where the last row is the Euler sequence of $\P(\Omega_Y)$, and
    take $D:=d\rho^{-1}(L)$. Equivalently, we may define $D\subset
    T_X$ by the rule $D_x=d\rho_x^{-1}(\langle x\rangle)$, where
    $\langle x\rangle\subset T_{X,\rho(x)}$ denotes the vector
    subspace determined by $x$.

    Since the statement is local, we will assume that $X=Y\times C$
    and that the zeroes of $\lambda$ in $\P(\Omega_Y)$ can be written as
    the classes of elements $V=V_1+tV_2+t^2V_3\in T_Y$, $t\in \C$,
    where $V_1,V_2$ and $V_3$ are independent vector fields on $Y$. We
    will denote by $V, V_1,V_2$ and $V_3$ the corresponding vector
    fields in $X$ obtained by the identification $T_X\iso
    p_1^*T_Y\oplus p_2^*T_C$, where $p_1$($=\rho$) and $p_2$ are the two natural projections.  Let us fix the following notation:
    $$
    [V_i,V_j]=\sum_{i=1}^3f_{ij}^kV_k.
    $$
    With the above notation, trivial computations show that $D$ is
    generated $\partial/\partial t$ and $V$, $\cE=[D,D]$ is a rank $3$
    vector bundle generated by $D$ and $V_2+2tV_3$, defined by the
    $1$-form
    $$\theta=t^2\omega_1-2t\omega_2+\omega_3,$$
    where the $\omega_i$'s
    denote $1$-forms dual to the $V_i$'s at each point.

    A straightforward computations shows that $\theta\wedge d\theta$
    takes the form:
    $$
    \theta\wedge d\theta=a\omega_1\wedge\omega_2\wedge\omega_3+bdt
    \wedge\omega_1\wedge\omega_2+cdt\wedge\omega_1\wedge
    \omega_3-2dt\wedge\omega_2\wedge\omega_3,
    $$
%%     $$
%%     \theta\wedge
%%     d\theta=\alpha\omega_1\wedge\omega_2\wedge\omega_3-2t^2dt
%%     \wedge\omega_1\wedge\omega_2+2tdt\wedge\omega_1\wedge
%%     \omega_3-2dt\wedge\omega_2\wedge\omega_3
%%     $$
%%     where
%%     $$
%%     \alpha=-t^4f_{23}^1+2t^3(f_{23}^2-f_{13}^1)-t^2(f_{23}^3+
%%     f_{12}^1-4f_{13}^2)+2t(f_{12}^2-f_{13}^3)+f_{12}^3$$
%%     Furthermore, the kernel $L$ of $\cE$ is generated by
%%     $$
%%     f_{12}^3+(2f_{13}^3-f_{12}^2)t+(f_{23}^3-2f_{13}^2)t^2-f_{23}^2t^3)
%%     \partial/\partial t-2V.
%%     $$
    and in particular it is everywhere non-zero.

    Easy examples of varieties $Y$ verifying the above property are
    the $3$-dimensional quadric $Q\in\P^4$ and abelian varieties. In
    the first case, the variety $X$ coincides with the flag manifold
    of pairs point-line contained in the quadric $Q$, and this is
    nothing but the projectivization of a contact structure on $\P^3$
    constructed as in Example \ref{ex:main}. However, the second case
    is of completely different nature, since $X$ does not admit a
    second $\P^1$ bundle structure (see also Proposition \ref{prop:quadric}).}
\end{ex}

Note that Jahnke and Radloff have recently obtained the complete list of holomorphic conformal structures on complex $3$-folds (cf. \cite{JR}), so providing the complete list of Engel structures of this type.
\begin{thm} \label{thm:conformal}
  The list of projective threefolds with a holomorphic conformal
  structure is as follows:
\begin{enumerate}
\item $Q_3$;
\item \'etale quotients of abelian threefolds;
\item threefolds with universal covering space the $3$-dimensional
  Lie ball $\D^{IV}_3$.
\end{enumerate}
\end{thm}
We denote here by $\D^{IV}_3$ the bounded symmetric domain dual to the
$3$-dimensional hyperquadric $Q_3$.

%%%%%%%%%%%%%%%%%%%%%%%%%%%%%%%%%%%%%%%%%%%%%%%%%%%%%%%%%%%%%

%%%%%%%%%%%%%%%%%%%%%%%%%%%%%%%%%%%%%%%%%%%%%%%%%%%%%%%%%%%%%
\section{Uniruledness of Engel manifolds}\label{sec:thm}

It is well known that the existence of rational curves on a variety $X$ depends on the positivity of its anticanonical divisor. By analogy with the contact case (see \cite{De}) we begin by inferring positivity properties of the anticanonical divisor of an Engel manifold from the non-integrability of the structure. From that we show that under certain positivity conditions on $L$ or $D/L$, the examples described in the previous section are the only possible Engel structures.

\subsection{The canonical class of an Engel manifold}\label{subsec:canon}

%\begin{prop}
%  Let $X$ be a projective variety admitting an Engel structure, then
%  the canonical divisor is not nef.
%\end{prop}
%
%\begin{proof}
We begin by recalling the following theorem by J.P. Demailly (cf.
\cite{De}):
\begin{thm}\label{thm:Dem}
  Let $X$ be a K\"ahler manifold. Assume that there is a
  pseudo-effective line bundle $N$ on $X$ and a nonzero holomorphic
  section $\theta \in H^0(X, \Omega_X\otimes N^{-1})$. Then the subsheaf
  defined by the kernel of $\theta$ defines a holomorphic foliation of
  codimension $1$ in $X$, that is, $\theta \wedge d\theta=0$.
\end{thm}

We will also make use of the characterization of pseudo-effective
divisors in terms of movable curves. A curve $C$ in $X$ is called {\it
  movable} if there exists an irreducible algebraic family of curves
containing $C$ and dominating $X$. Boucksom, Demailly, Peternell and
Paun have recently shown (cf. \cite{BDPP}) that a divisor is
pseudo-effective if and only if it has non negative degree on every movable
curve.

Applying this result to the section $\theta\in H^0(X,\Omega_X\otimes
L')$ defining $\cE$ for an Engel structure we obtain the following straightforward result:

\begin{lem}\label{lem:unirul}
  Let $X$ be a complex manifold admitting a holomorphic Engel
  structure $L\subset D\subset \cE$. Then either $L^{-1}$ or $(D/L)^{-1}$ are not pseudo-effective. Moreover, assume that any of the following
  properties is fulfilled:
\begin{itemize}
\item $L$ is pseudo-effective,
\item $\det(D)$ is pseudo-effective.
\end{itemize}
Then the canonical divisor of $X$ is not pseudo-effective and,
equivalently, $X$ is uniruled.
\end{lem}
\begin{proof}
The first assertion follows directly from Theorem \ref{thm:Dem}, which in our case tells us that $(L')^{-1}=\big(L\otimes(D/L)^2\big)^{-1}$ is not pseudo-effective.

In order to get the second, note that by Lemma \ref{lem:canon} we get $\omega_X^{-1}\cong L\otimes(L')^2\cong L'\otimes \det(D)$. Thus if
  $L$ or $\det(D)$ are pseudo-effective and $L'$ has positive degree
  on a movable curve $C$, then $\omega_X\cdot C<0$, allowing us to conclude
  using \cite{BDPP}. But the condition on $L'$ is equivalent (again by
  \cite{BDPP}) to the non pseudo-effectivity of $(L')^{-1}$, that we
  obtain from \ref{thm:Dem}.
\end{proof}
%
%\begin{prb}\label{prb:unirul}
%{\rm It would be nice to understand Demailly's theorem at the level of movable curves, i.e. to construct geometrically a movable curve of negative $N$-degree associated to a non-integrable section $\theta$. We believe that this would bring some light on the complete classification of Engel structures.}
%\end{prb}
%In our case the associated even contact structure defines a section
%$\theta\in H^0(\Omega_X \otimes L')$. Recall that $L'=L\otimes (D/L)^2$. Assume
%that $X$ is not an element of the two families of examples already classified,
%therefore $L^{-1}$ and $(D/L)^{-1}$ are nef. So
%$(L')^{-1}=L^{-1}\otimes (D/L)^{-1}$ is nef. Moreover, nef implies
%pseudo-effective \cite{BDPP} and so $(L')^{-1}$ is pseudo-effective. We are
%in the hypothesis of the previous Theorem and so $\theta \wedge
%d\theta=0$ which is a contradiction. Therefore the only possible cases
%are the two ones already studied.
%
%\end{proof}

\subsection{Case I: $L^{-1}$ is not pseudoeffective}

In this section we will show that the Engel structures described in Example \ref{ex:main} may be characterized by the non pseudo-effectivity of $L^{-1}$. Note that whenever $L$ is not numerically trivial, this property is weaker than the pseudo-effectivity of $L$ needed in Lemma \ref{lem:unirul}.

\begin{prop}\label{prop:P1struct}
  Let $X$ be a smooth complex projective variety and $\cL\subset T_X$
  a line subbundle of the tangent bundle of $X$. Assume that
  $\cL^{-1}$ is not nef. Then $X$ is isomorphic to a $\P^1$-bundle
  over a nonsingular variety $Z$.

%%   Let $\cE$ be an even contact strucure on the $4$-dimensional
%%   manifold $X$, with kernel $L$. Assume that $L^{-1}$ is not nef. Then
%%   there exists a $3$-dimensional manifold $Z$ and a rank $2$ vector
%%   bundle $E$ on $Z$ such that $X\iso\P(E)$.
\end{prop}

\begin{proof}
  The line subbundle $\cL\subset T_X$ defines a
  $1$-dimensional regular foliation on $X$. By assumption there exists
  a curve $C\subset X$ such that $\cL|_C$ is ample. It follows from
  \cite[Thm.~1]{KST} that the leaves of $\cL$ passing by points of $C$
  are algebraic and the general one is rationally connected, hence a
  $\P^1$.

  Now we claim that the general leaf of $\cL$ is isomorphic to $\P^1$.
  Fix a leaf $R$ passing by a general point of $C$. Since $R\iso\P^1$
  the holonomy of the foliation $\cL$ along $R$ is trivial and Reeb's
  Stability Theorem \cite[Thm.~IV.3]{CLN} provides a fundamental
  system of analytic neighborhoods of $R$ in $X$ that are saturated
  with respect to $\cL$. It follows that the neighboring leaves of
  $\cL$ are $\P^1$'s. This defines an analytic family of $\P^1$'s
  parametrized by a submanifold $S_0\subset\RatCurves^n(X)$. Now,
  every element in the Zariski closure inside
  $S\subset\RatCurves^n(X)$ parametrizes a rational curve tangent to
  $\cL$. Since the family parametrized by $S$ dominates $X$, we
  conclude the claim.

  Now we prove that every leaf of $\cL$ is a $\P^1$. Let $Z$ be
  the normalization of the closure of the family $S$ in $\Chow(X)$.
  Every element $z\in Z$ determines a rational cycle $C_z$ in $X$.
  Since $C_z$ is a limit of elements of $S$, it must be tangent to
  $\cL$ at a smooth point. But $L$ is regular, hence every cycle
  $C_z$ must be smooth, reduced and irreducible, and $Z$ is smooth at
  $z$. 

  Finally, since $L\cdot \pi^{-1}(z)=2$ for every $z\in Z$, Lemma \ref{lem:canon} tells us that $(D/L)^{-1}$ has degree $1$ on every fiber. This concludes the proof. 

\end{proof}

\begin{prop}\label{prop:cartan}
  Let $X$ be a smooth projective $2n$-fold admitting an even contact
  structure $\cE\subset T_X$ with kernel $L$. If $L^{-1}$ is not nef (in particular, if $L^{-1}$ is not pseudoeffective), then
  $X$ is a $\P^1$-bundle over a contact $(2n-1)$-fold $Z$.
  Furthermore, if $n=4$ and $\cE$ is an Engel structure, then $X$ is the
  projectivization of a contact distribution on $Z$.
\end{prop}

\begin{proof}
  Consider the variety $Z$ and the morphism $\pi:X\to Z$ constructed in
  Proposition \ref{prop:P1struct} by integration of the foliation
  $L\subset T_X$.
  Consider the exact sequence
  $$0\rightarrow\cE/L\longrightarrow T_X/L\cong
  \pi^*T_Z\longrightarrow L'\rightarrow 0.$$

  Arguing as in Prop.~\ref{prop:fibration}, $L'$ and $\cE/L$ are
  trivial on every fiber of $\pi$. By \cite[III, 12.9]{Hart2} it follows that $\cE/L=\pi^*E$,
  $L'=\pi^*L'_Z$, where $E$ and $L'_Z$ are locally free sheaves of rank
  $2n-2$ and $1$, respectively. Applying $\pi_*$ to the exact sequence
  above and considering that $R^1\pi_*\cE/L=0$, we obtain an exact
  sequence
  $$0\rightarrow E\longrightarrow
  T_Z\stackrel{\varphi}{\longrightarrow} L'_Z\rightarrow 0.$$
  But $\varphi$ defines a contact structure, because pulling it back and composing it with $d\pi$ we obtain the original even contact structure on $X$ (see Example \ref{ex:main}).

  Finally if $\cE$ is associated with an Engel structure $D$, then $D/L$ and $L\otimes D/L$ have degree
  $-1$ and $1$, respectively, on the fibers of $\pi$. It follows that
  $\pi_*D/L=R^1\pi_*D/L=0$ and pushing down the exact sequence:
  $$
  0\rightarrow D/L\longrightarrow \cE/L=\pi^*E\longrightarrow
  L\otimes D/L\rightarrow 0
  $$
  we get $E\cong \pi_*(L\otimes D/L)$, therefore $X\cong\P(E)$.
\end{proof}

\subsection{Case II: $(D/L)^{-1}$ is not pseudoeffective}

We now study the case of $(D/L)^{-1}$ being not pseudoeffective. In this case we need to make a second assumption in order to characterize Lorentzian tubes.

\begin{prop}\label{prop:tubes}
  Let $X$ be a smooth projective $4$-fold admitting an Engel structure
  $L\subset D\subset\cE$ such that $(D/L)^{-1}$ is not nef.  Assume
  that $L$ is a direct summand of $D$. Then there exist a smooth
  $3$-fold $Y$ and a holomorphic conformal structure $\lambda\in
  H^0(Y, S^2\Omega_Y\otimes M)$, whose set of zeroes is isomorphic to
  $X$.
\end{prop}

\begin{proof}
  The assumption on $L$ allows us to consider $D/L$ as a line
subbundle of $T_X$. Applying Prop. \ref{prop:P1struct} we deduce that $X$ is a $\P^1$-bundle over some smooth variety $Y$. Denote by
  $\rho$ the natural projection from $X$ to $Y$. Since the relative
  tangent bundle to $\rho$ is $D/L$, there is an inclusion $L\to
  \pi^*T_Y$, which defines an inclusion
  $X\hookrightarrow\P(\Omega_Y)$. But then $D/L$ has degree $2$ on every fiber $C$ of $\rho$, hence Lemma \ref{lem:canon} tells us that $L\cdot C=-2$. This concludes the proof.
\end{proof}

Finally we show that the two classes of varieties admitting an Engel
structure previously described have a unique common element, namely the flag manifold
$F_Q(0,1)$, where $Q$ is a quadric threefold. That is the universal
family of lines contained in $Q$. Note that this variety admits two
$\P^1$-bundle structures
$$
\xymatrix{F_Q(0,1)\ar[r]^{\rho} \ar[d]^{\pi}&Q\\\P^3&}
$$

\begin{prop}\label{prop:quadric}
  Let $X$ be a smooth projective manifold admitting and Engel
  structure $L\subset D\subset \cE$. Assume that $L^{-1}$ and
  $(D/L)^{-1}$ are not nef, and that $D=L\oplus D/L$. Then $X$ is isomorphic to the flag
  manifold $F_Q(0,1)$.
\end{prop}

\begin{proof}
  Consider the two $\P^1$-bundle structures $\pi:X\to Z$ and
  $\rho:X\to Y$ constructed in the previous propositions We will make
  use of the notation of the proofs appearing there. If we prove that
  $Z=\P^3$, then $X$ is the projectivization of a contact structure on
  $Z$. But it is classically known that then $X\cong F_Q(0,1)$. 
  
  In order to prove that $Z=\P^3$ it is enough to check that $X$ contains a family of rational curves $\cM=\{C_t\}$, with $\cM$ proper and splitting type $$T_Z|_{C_t}=\cO_{\P^1}(2)\oplus\cO_{\P^1}(1)^{\oplus 2}\mbox{ for all }C_t.$$ 
  If that property is fulfilled, then every two points of $Z$ might be joined by a curve of the family $\cM$ (cf. \cite[Prop. 4.8]{De}), and in particular \cite[Lemma 1]{O} tells us that $Z$ is a Fano threefold of Picard number $1$. Moreover, $-K_Z\cdot C_t=4$ implies that the index of $Z$ is $4$, hence $Z\cong\P^3$. 
  
  Consider a fiber $C:=\rho^{-1}(z)$ for
  any $z\in Z$. By the computations we have done in Examples \ref{ex:main} and \ref{ex:second}, it suffices to show that
  %$\pi^*T_Z|_C\cong\cO(2)\oplus\cO(1)^{\oplus 2}$, or equivalently that
  $\big(T_X/D\big)|_C\cong\cO(1)^{\oplus 2}$. But $T_X/D\cong
  \big(T_{\P(\Omega_Y)|Y}(-1)\big)|_X$ and this is isomorphic to $\cO(1)^{\oplus 2}$
  because $T_{\P(\Omega_Y),Y}|_C\cong T_{\P^2,C}$ and $C\subset\P^2$
  is a conic.
\end{proof}

%%%%%%%%%%%%%%%%%%%%%%%%%%%%%%%%%%%%%%%%%%%%%%%%%%%%%%%%%%%%%

%%%%%%%%%%%%%%%%%%%%%%%%%%%%%%%%%%%%%%%%%%%%%%%%%%%%%%%%%%%%%

\hrule \medskip
\par\noindent

Addresses:

\par\noindent{\tt e-mail fpresas@imaff.cfmac.csic.es}
\par\noindent{\tt e-mail luis.sola@urjc.es}

\end{document}